\documentclass[a4paper,11pt]{amsart}

\usepackage[top=30truemm,bottom=25truemm,left=35truemm,right=35truemm]{geometry}

\usepackage{amsmath}
\usepackage{amssymb}
\usepackage{amsthm}
\usepackage{graphics}
\usepackage[abbrev,alphabetic]{amsrefs}
\usepackage{amscd}
\usepackage[dvipdfmx]{graphicx}
\usepackage{color}
\usepackage{ulem}
\usepackage{enumitem}

\usepackage[all,cmtip]{xy}
\usepackage{xypic}

\title[Precise inversion of adjunction in dimension four]
{Precise inversion of adjunction in dimension four}

\author{Yusuke Nakamura}
\address{Graduate School of Mathematics, Nagoya University, Furo-cho, Chikusa-ku, Nagoya, 464-8602, Japan.}
\email{y.nakamura@math.nagoya-u.ac.jp}
\urladdr{https://sites.google.com/site/ynakamuraagmath/}

\author{Kohsuke Shibata}
\address{Faculty of Science and Technology, Tokushima University, 2-1 Minamijosanjima,Tokushima, 770-8506, Japan.}
\email{shibata.kohsuke@tokushima-u.ac.jp}

\subjclass[2020]{Primary 14E30; Secondary 14B05, 14N30}

\keywords{minimal log discrepancy, precise inversion of adjunction, PIA conjecture}

\newtheorem{thm}{Theorem}[section]
\newtheorem{lem}[thm]{Lemma}

\newtheorem{prop}[thm]{Proposition}

\theoremstyle{definition}
\newtheorem{defi}[thm]{Definition}

\newtheorem{conj}[thm]{Conjecture}

\theoremstyle{remark}
\newtheorem{rmk}[thm]{Remark}

\newtheorem*{ackn}{Acknowledgements}
\newtheorem*{aidisc}{Declaration of generative AI use}

\begin{document}

\maketitle

\begin{abstract}
We prove the precise inversion of adjunction conjecture for minimal log discrepancies on four-dimensional varieties.
\end{abstract}

\section{Introduction}

Throughout this paper, we work over an algebraically closed field $k$ of characteristic zero.

Minimal log discrepancies are fundamental invariants of singularities in birational geometry.
They measure the complexity of singularities and play an important role in the minimal model program, particularly in the context of the termination of flips (cf.\ \cite{Sho04}).
A basic tool for studying singularities by induction on dimension is the inversion of adjunction.

\begin{thm}[Inversion of adjunction, {\cite{92}*{17.6}, \cite{Kaw07}}]\label{thm:IA}
Let $(X, \Delta)$ be a log pair over $k$ and let $S$ be a normal Cartier prime divisor on $X$ such that $S \not\subset \operatorname{Supp}(\Delta)$.
\begin{enumerate}
\item $(X, S + \Delta)$ is plt near $S$ if and only if $(S, \Delta|_S)$ is klt.
\item $(X, S + \Delta)$ is log canonical near $S$ if and only if $(S, \Delta|_S)$ is log canonical.
\end{enumerate}
\end{thm}

The precise inversion of adjunction (PIA) conjecture asserts that such equivalences hold at the level of minimal log discrepancies.

\begin{conj}[PIA conjecture, {cf.\ \cite{92}*{17.3.1}}]\label{conj:PIA}
Let $(X, \Delta)$ be a log pair over $k$ and let $S$ be a normal Cartier prime divisor on $X$ such that $S \not \subset \operatorname{Supp} (\Delta)$.
Then for every closed point $x \in S$,
\[
\operatorname{mld}_x(X, S + \Delta) = \operatorname{mld}_x(S, \Delta|_S)
\]
holds.
\end{conj}

The inequality ``$\le$'' always holds (\cite{92}*{17.2}), and the opposite inequality is the essential part of the conjecture.
Conjecture~\ref{conj:PIA} is known to hold when $\dim X = 2$ \cite{Sho93} and when $\dim X = 3$ \cite{NS26b}*{Proposition 6.1}, when $X$ is smooth \cite{EMY03}, and when $X$ is a local complete intersection \cite{EM04}.
Moreover, it is known to hold when $X$ is a quotient of a complete intersection singularity and $S$ is klt at $x$ \cite{NS22, NS25, NS26a}.

While Theorem~\ref{thm:IA} requires no assumption on the singularities of $S$, the PIA conjecture fails in general without such an assumption: in \cite{NS26b}, a counterexample is constructed for every $\dim X \ge 5$, in which $S$ is not klt at $x$.
Thus the conjecture remained open only in dimension four (\cite{NS26b}*{Question 6.2}).

The purpose of this paper is to settle this remaining case.

\begin{thm}[$=$ Theorem~\ref{thm:PIA_dim4}]\label{thm:PIA_dim4_intro}
Conjecture~\ref{conj:PIA} holds when $\dim X = 4$.
\end{thm}

We prove Theorem~\ref{thm:PIA_dim4_intro} by dividing it into the following two cases:
\begin{itemize}
\item $\operatorname{mld}_x(X, S + \Delta) < 1$, and
\item $\operatorname{mld}_x(X, S + \Delta) \ge 1$.
\end{itemize}

In the first case, the conjecture can be proved in arbitrary dimension by using Moraga's result \cite{Mor20} on the extraction of a divisor with log discrepancy in $[0, 1)$ over a log canonical pair (see Theorem~\ref{thm:Moraga}).

\begin{prop}[$=$ Proposition~\ref{prop:PIA_low}]\label{prop:PIA_low_intro}
Conjecture~\ref{conj:PIA} holds when $\operatorname{mld}_x(X, S + \Delta) < 1$.
\end{prop}

In the second case, we first reduce to the case where $X$ is $\mathbb{Q}$-Gorenstein by taking the canonical modification of $X$ (Proposition~\ref{prop:reduction_QGor}).
We may then assume that $\operatorname{mld}_x(S, \Delta |_S) > 1$.
In this case, $S$ is klt at $x$, and the divisor $S'$ on the index-one cover $X'$ of $X$ at $x$ corresponding to $S$ has a Gorenstein canonical singularity with minimal log discrepancy at least two.
By the classification of three-dimensional Gorenstein canonical singularities due to Reid \cite{Rei80}, $S'$ is either smooth or has a cDV singularity, and it follows that $X'$ is smooth or has a hypersurface singularity.
We are thus reduced to the setting of the precise inversion of adjunction for quotients of complete intersection singularities, established in \cite{NS22, NS25, NS26a}.
We use it in the following form.

\begin{thm}[$=$ Theorem~\ref{thm:quotient_PIA}]\label{thm:quotient_PIA_intro}
Let $X'$ be a normal quasi-projective variety, let $G$ be a finite group acting on $X'$ freely in codimension one, and let $x' \in X'$ be a closed point such that $X'$ is a local complete intersection at $x'$.
Let $q \colon X' \to X := X'/G$ be the quotient morphism.
Let $(X, \Delta)$ and $S$ be as in Conjecture~\ref{conj:PIA} with $x := q(x') \in S$, and suppose that $S$ is klt at $x$.
Then Conjecture~\ref{conj:PIA} holds for $(X, \Delta)$, $S$ and $x$.
\end{thm}

Theorem~\ref{thm:quotient_PIA_intro} is an extension of \cite{NS26a}*{Theorem 7.3}.
The latter assumes that $X'$ is embedded $G$-equivariantly into an affine space with a linear $G$-action and is defined there as a complete intersection by $G$-semi-invariant equations, whereas Theorem~\ref{thm:quotient_PIA_intro} only requires $X'$ to be a local complete intersection at $x'$ (see Remark~\ref{rmk:comparison_NS} for details).
When $G$ is cyclic, which is the case needed for Theorem~\ref{thm:PIA_dim4_intro}, we can construct a $G$-equivariant embedding of $X'$ near $x'$ into an affine space with a linear $G$-action in which $X'$ is defined by a regular sequence of $G$-semi-invariant polynomials (Lemma~\ref{lem:equiv_emb}).
This allows us to apply \cite{NS26a}*{Theorem 7.3}.
The case of a general finite group $G$ is then reduced to the cyclic case by \cite{NS22}*{Remark 4.13}.

As a byproduct of the proof in the second case, we obtain the following structure theorem for four-dimensional singularities $x \in X$ admitting a normal Cartier divisor $S$ with $\operatorname{mld}_x(X, S) > 1$.

\begin{thm}[$=$ Theorem~\ref{thm:hyperquotient}]\label{thm:hyperquotient_intro}
Let $X$ be a normal $\mathbb{Q}$-Gorenstein variety with $\dim X = 4$, let $S$ be a normal prime Cartier divisor on $X$, and let $x \in S$ be a closed point.
Suppose that $\operatorname{mld}_x(X, S) > 1$.
Then $X$ has a hyperquotient singularity at $x$. More precisely, the index-one cover of $X$ at $x$ is smooth or has a hypersurface singularity at the point over $x$.
\end{thm}

\medskip
This paper is organized as follows.
In Section~\ref{section:prelim}, we recall the definitions of minimal log discrepancies for log pairs and for $\mathbb{R}$-ideal sheaves.
In Section~\ref{section:arbitrary_dim}, we prove two results valid in arbitrary dimension: Conjecture~\ref{conj:PIA} in the case $\operatorname{mld}_x(X, S + \Delta) < 1$ (Proposition~\ref{prop:PIA_low}), and the reduction of Conjecture~\ref{conj:PIA} to the case where $X$ is $\mathbb{Q}$-Gorenstein (Proposition~\ref{prop:reduction_QGor}).
In Section~\ref{section:quotient}, we construct equivariant embeddings of local complete intersection singularities with cyclic group actions (Lemma~\ref{lem:equiv_emb}) and prove Theorem~\ref{thm:quotient_PIA}.
In Section~\ref{section:dim4}, we treat the $\mathbb{Q}$-Gorenstein case with $\operatorname{mld}_x(X, S + \Delta) \ge 1$ in dimension four (Proposition~\ref{prop:PIA_QGor_dim4}), and complete the proofs of Theorems~\ref{thm:PIA_dim4} and \ref{thm:hyperquotient}.

\begin{ackn} 
The first author is partially supported by Inamori Foundation and by JSPS KAKENHI No.\ 22K13888. 
The second author is partially supported by JSPS KAKENHI No.\ 23K12958. 
\end{ackn}

\begin{aidisc}
This project started around July 2026, and its theoretical part was completed by mid-September 2026.
The project originated from the authors' idea, which we consider to be the main idea of this paper: to divide the problem into the cases $\operatorname{mld}_x(X, S + \Delta) < 1$ and $\operatorname{mld}_x(X, S + \Delta) \ge 1$, and, in the latter case, to show that $X$ has a hyperquotient singularity by the classification of the singularities of $S$, so that the authors' previous results \cite{NS22, NS25, NS26a} can be applied.
Apart from this main idea, about half of the technical parts of the proofs were developed through interactive discussions with generative AI tools, namely ChatGPT and Claude. 
In particular, the observation that \cite{NS26a}*{Theorem 7.3} can be extended to Theorem~\ref{thm:quotient_PIA} relied substantially on these discussions.
The proof of Lemma~\ref{lem:equiv_emb}, which is needed for this extension, also relied substantially on them.
We also used ChatGPT and Claude to check the proofs. 
The organization of the paper and the exposition of the proofs were written by the authors, and Claude was then used for English proofreading. 
The authors have verified all the arguments in this paper and take full responsibility for its content.
\end{aidisc}

\section{Preliminaries}\label{section:prelim}

\subsection{Notation}

\begin{itemize}
\item We basically follow the notation and terminology in \cite{KM98} and \cite{Kol13}.

\item Throughout this paper, $k$ is an algebraically closed field of characteristic zero. A \textit{variety} is an integral scheme that is separated and of finite type over $k$.

\item A variety $X$ is a \textit{local complete intersection at} a closed point $x$ if $\mathcal{O}_{X, x}$ is a complete intersection ring. 

\item A variety $X$ is said to have a \textit{hyperquotient singularity} at a closed point $x$ if, in a neighborhood of $x$, the variety $X$ is the quotient of a variety which is smooth or has a hypersurface singularity at a point over $x$ by a finite group.

\item Let $G$ be a finite group acting on a normal variety $V$. For $g \in G$, we denote by $V^g$ the fixed locus of $g$. We say that $G$ acts on $V$ \textit{freely in codimension one} if $\operatorname{codim}_V V^g \ge 2$ for every $1 \neq g \in G$. 
This holds if and only if $G$ acts faithfully on $V$ and the quotient morphism $V \to V/G$ is \'etale in codimension one. 

\item Let $G$ be a finite group acting on a $k$-algebra $A$ by $k$-algebra automorphisms. An element $a \in A$ is said to be \textit{$G$-semi-invariant} if there exists a character $\chi \colon G \to k^*$ such that $g(a) = \chi(g) a$ for every $g \in G$, and \textit{$G$-invariant} if moreover $\chi$ is trivial.
\end{itemize}

\subsection{Minimal log discrepancies}

A \textit{log pair} $(X, \Delta)$ over $k$ consists of a normal variety $X$ over $k$ and an effective $\mathbb{R}$-divisor $\Delta$ on $X$ such that $K_X + \Delta$ is $\mathbb{R}$-Cartier. 

\begin{defi}\label{defi:mld}
Let $(X, \Delta)$ be a log pair over $k$.
\begin{enumerate}
\item
Let $f \colon Y \to X$ be a proper birational morphism from a normal variety $Y$, and let $E$ be a prime divisor on $Y$.
Let $\Delta_Y$ be the $\mathbb{R}$-divisor on $Y$ defined by $K_Y + \Delta_Y = f^*(K_X + \Delta)$ and $f_*K_Y = K_X$.
The \textit{log discrepancy} of $(X, \Delta)$ at $E$ is defined as
\[
a_E(X, \Delta) := 1 - \operatorname{coeff}_E \Delta_Y.
\]
The image $f(E)$ is called the \textit{center} of $E$ on $X$ and is denoted by $c_X(E)$.

\item
Let $\eta$ be a scheme-theoretic point of $X$ with $\operatorname{codim} \eta \ge 1$, and set
\[
\alpha_\eta := \inf_{c_X(E) = \overline{\{\eta\}}} a_E(X, \Delta),
\]
where the infimum is taken over all prime divisors $E$ over $X$ with $c_X(E) = \overline{\{\eta\}}$.
The \textit{minimal log discrepancy} of $(X, \Delta)$ at $\eta$ is defined by
\[
\operatorname{mld}_\eta(X, \Delta) :=
\begin{cases}
\alpha_\eta & \text{if } \operatorname{codim} \eta \ge 2, \\
\alpha_\eta & \text{if } \operatorname{codim} \eta = 1 \text{ and } \alpha_\eta \ge 0, \\
-\infty & \text{if } \operatorname{codim} \eta = 1 \text{ and } \alpha_\eta < 0.
\end{cases}
\]
For a closed point $x \in X$, we write $\operatorname{mld}_x(X, \Delta)$ for the minimal log discrepancy at $x$.

\item
For a closed subset $W \subset X$, we define
\[
\operatorname{mld}_W(X, \Delta) := \inf_{\eta} \operatorname{mld}_\eta(X, \Delta),
\]
where the infimum is taken over all scheme-theoretic points $\eta \in W$ with $\operatorname{codim} \eta \ge 1$.

\item
We say that a prime divisor $E$ over $X$ \textit{computes} $\operatorname{mld}_\eta(X, \Delta)$ if $c_X(E) = \overline{\{\eta\}}$ and
\[
\begin{cases}
a_E(X, \Delta) = \operatorname{mld}_\eta(X, \Delta) & \text{if } \operatorname{mld}_\eta(X, \Delta) \ge 0, \\[0.5em]
a_E(X, \Delta) < 0 & \text{if } \operatorname{mld}_\eta(X, \Delta) = -\infty.
\end{cases}
\]
\end{enumerate}
\end{defi}

We also use minimal log discrepancies with respect to $\mathbb{R}$-ideal sheaves.
An \textit{$\mathbb{R}$-ideal sheaf} on a variety $X$ is a formal product $\mathfrak{a} = \prod_{i=1}^s \mathfrak{a}_i^{\tau_i}$, where $\mathfrak{a}_1, \dots, \mathfrak{a}_s$ are nonzero coherent ideal sheaves on $X$ and $\tau_1, \dots, \tau_s$ are positive real numbers.
For a morphism $g \colon Y \to X$, we set $\mathfrak{a}\mathcal{O}_Y := \prod_{i=1}^s (\mathfrak{a}_i \mathcal{O}_Y)^{\tau_i}$, and for a prime divisor $E$ over $X$, we set $\operatorname{ord}_E(\mathfrak{a}) := \sum_{i=1}^s \tau_i \operatorname{ord}_E(\mathfrak{a}_i)$.

\begin{defi}\label{defi:mld_ideal}
Let $X$ be a normal $\mathbb{Q}$-Gorenstein variety over $k$ and let $\mathfrak{a}$ be an $\mathbb{R}$-ideal sheaf on $X$.
For a prime divisor $E$ over $X$, the \textit{log discrepancy} of $(X, \mathfrak{a})$ at $E$ is defined as
\[
a_E(X, \mathfrak{a}) := a_E(X, 0) - \operatorname{ord}_E(\mathfrak{a}).
\]
The minimal log discrepancy $\operatorname{mld}_\eta(X, \mathfrak{a})$ at a scheme-theoretic point $\eta \in X$ with $\operatorname{codim} \eta \ge 1$ is defined as in Definition~\ref{defi:mld}~(2), with $a_E(X, \Delta)$ replaced by $a_E(X, \mathfrak{a})$.
\end{defi}

\begin{rmk}\label{rmk:ideal_setting}
Let $(X, \Delta)$ be a log pair such that $X$ is $\mathbb{Q}$-Gorenstein.
Then $\Delta$ is an effective $\mathbb{R}$-Cartier divisor, and hence it can be written as $\Delta = \sum_j \lambda_j D_j$ with $\lambda_j > 0$ and effective Cartier divisors $D_j$ (see \cite{Kol23}*{(11.45) and Proposition 11.47}).
We obtain an $\mathbb{R}$-ideal sheaf
$\mathfrak{a}_\Delta := \prod_j \mathcal{O}_X(-D_j)^{\lambda_j}$ on $X$ such that
$\operatorname{ord}_E(\mathfrak{a}_\Delta) = \operatorname{ord}_E(\Delta)$ for every prime divisor $E$ over $X$.
Thus we have $\operatorname{mld}_\eta(X, \Delta) = \operatorname{mld}_\eta(X, \mathfrak{a}_\Delta)$
for every scheme-theoretic point $\eta \in X$ with $\operatorname{codim} \eta \ge 1$.
\end{rmk}

\begin{rmk}\label{rmk:easy}
In the setting of Conjecture~\ref{conj:PIA}, the inequality
\[
\operatorname{mld}_x(X, S + \Delta) \le \operatorname{mld}_x(S, \Delta |_S)
\]
always holds (see \cite{92}*{17.2}).
Moreover, Conjecture~\ref{conj:PIA} holds if $\operatorname{mld}_x(X, S + \Delta) = -\infty$.
Indeed, in this case $(X, S + \Delta)$ is not log canonical at $x$, so $(S, \Delta |_S)$ is not log canonical at $x$ by inversion of adjunction for log canonicity (Theorem \ref{thm:IA} (2)), and hence $\operatorname{mld}_x(S, \Delta |_S) = -\infty$.
\end{rmk}

\section{Results in arbitrary dimension}\label{section:arbitrary_dim}

\subsection{The case of small minimal log discrepancy}

In this subsection, we prove Conjecture~\ref{conj:PIA} when $\operatorname{mld}_x(X, S + \Delta) < 1$, in arbitrary dimension.
We use the following result of Moraga, which allows us to extract a single divisor with log discrepancy less than one over a log canonical pair.

\begin{thm}[{\cite{Mor20}*{Theorem 1}}]\label{thm:Moraga}
Let $(X, \Delta)$ be a log canonical pair, where $X$ is quasi-projective, and let $E$ be an exceptional prime divisor over $X$ with $a_E(X, \Delta) < 1$.
Then there exists a projective birational morphism $f \colon Y \to X$ from a normal variety $Y$ whose exceptional locus is $E$, such that $E$ is $\mathbb{Q}$-Cartier.
\end{thm}

\begin{proof}
By \cite{Kol13}*{Proposition 2.21}, we may assume that $\Delta$ is a $\mathbb{Q}$-divisor.  
Applying \cite{Mor20}*{Theorem 1} to $\mathcal{V} = \{\operatorname{ord}_E\}$, we obtain a projective birational morphism $f \colon Y \to X$ with $\operatorname{Exc}(f) = E$ such that $E$ is $\mathbb{Q}$-Cartier. 
\end{proof}

The following proposition was proved in \cite{Kaw21}*{Theorem 7.5} when $X$ is klt.
By using Theorem~\ref{thm:Moraga} in place of the extraction for klt pairs, the same argument works in the log canonical case.

\begin{prop}\label{prop:PIA_low}
Conjecture~\ref{conj:PIA} holds when $\operatorname{mld}_x(X, S + \Delta) < 1$.
\end{prop}

\begin{proof}
Let $(X, \Delta)$, $S$ and $x$ be as in Conjecture~\ref{conj:PIA}, and suppose that $\operatorname{mld}_x(X, S + \Delta) < 1$.
By Remark~\ref{rmk:easy}, it suffices to show
\[
\operatorname{mld}_x(S, \Delta |_S) \le \operatorname{mld}_x(X, S + \Delta).
\]
Furthermore, we may assume, after shrinking $X$ to an affine neighborhood of $x$, that $(X, S + \Delta)$ is log canonical. 
Let $E$ be a prime divisor over $X$ computing $\operatorname{mld}_x(X, S + \Delta)$ (see Definition \ref{defi:mld}(4)). 
Put
\[
a := a_E(X, S + \Delta) = \operatorname{mld}_x(X, S + \Delta),
\]
so that $0 \le a < 1$.
Let $f \colon Y \to X$ be the extraction of $E$ given by Theorem~\ref{thm:Moraga}, and let $S_Y$ and $\Delta_Y$ be the strict transforms of $S$ and $\Delta$. Then
\[\tag{$\heartsuit$}
K_Y + S_Y + \Delta_Y + (1 - a) E = f^*(K_X + S + \Delta),
\]
and the pair $(Y, S_Y + \Delta_Y + (1-a)E)$ is log canonical.

Since $f^{-1}(x) = E$ and $x \in S$, we have $S_Y \cap E \neq \emptyset$.
Since $E$ is $\mathbb{Q}$-Cartier, every irreducible component of $S_Y \cap E$ is a prime divisor on $S_Y$ mapped to $x$.
Let $Z$ be one of them, and let $\eta_Z$ denote its generic point.

Since $(Y, S_Y + \Delta_Y + (1-a)E)$ is log canonical and $\eta_Z \in \operatorname{Supp}\bigl(\Delta_Y + (1-a)E\bigr)$, we are in case (16.6.3) of \cite{92}*{16.6} (see the proof of \cite{92}*{16.7}).
In particular, $S_Y$ is normal at $\eta_Z$, and hence $Z$ defines a prime divisor over $S$ with center $x$.
Restricting ($\heartsuit$) to $S_Y$ by \cite{92}*{16.5}, we obtain
\[
K_{S_Y} + \operatorname{Diff}_{S_Y}\bigl(\Delta_Y + (1-a)E\bigr) = (f|_{S_Y})^*(K_S + \Delta |_S)
\]
near $\eta_Z$, and therefore
\[
\operatorname{mld}_x(S, \Delta |_S) \le a_Z(S, \Delta |_S) = 1 - \operatorname{coeff}_{\eta_Z} \operatorname{Diff}_{S_Y}\bigl(\Delta_Y + (1-a)E\bigr).
\]
Hence, in order to show $\operatorname{mld}_x(S, \Delta |_S) \le a$, it suffices to prove
\[
\operatorname{coeff}_{\eta_Z} \operatorname{Diff}_{S_Y}\bigl(\Delta_Y + (1-a)E\bigr) \ge 1 - a.
\]
By \cite{92}*{16.7}, we have
\[
\operatorname{coeff}_{\eta_Z} \operatorname{Diff}_{S_Y}\bigl(\Delta_Y + (1-a)E\bigr)
= 1 - \frac{1}{m} + \frac{r(1-a)}{m} + \sum_i \frac{r_i b_i}{m}
\]
for some integers $m, r \ge 1$ and $r_i \ge 0$, where $\Delta_Y = \sum_i b_i B_i$.
Hence
\[
\operatorname{coeff}_{\eta_Z} \operatorname{Diff}_{S_Y}\bigl(\Delta_Y + (1-a)E\bigr) \ge 1 - \frac{1}{m} + \frac{r(1-a)}{m} \ge 1 - a,
\]
which completes the proof.
\end{proof}

\subsection{Reduction to the $\mathbb{Q}$-Gorenstein case}

In this subsection, we show that Conjecture~\ref{conj:PIA} in a fixed dimension follows from the case where $X$ is $\mathbb{Q}$-Gorenstein (Proposition~\ref{prop:reduction_QGor}).
For this purpose, we use the canonical modification of $X$, which provides a $\mathbb{Q}$-Gorenstein model of $X$  (Lemma~\ref{lem:can_mod}).

\begin{lem}\label{lem:can_mod}
Let $(X, \Delta)$ be a log pair over $k$ and let $S$ be a normal Cartier prime divisor on $X$ such that $(X, S + \Delta)$ is log canonical.
Then there exist a projective birational morphism $f \colon Y \to X$ and an effective $\mathbb{R}$-divisor $\Gamma$ on $Y$ satisfying the following conditions, where $S_Y$ denotes the strict transform of $S$:
\begin{enumerate}
\item $Y$ is $\mathbb{Q}$-Gorenstein, 
\item $S_Y$ is a normal Cartier prime divisor and $f^{-1}(S) = S_Y$, and 
\item $K_Y + S_Y + \Gamma = f^*(K_X + S + \Delta)$.
\end{enumerate}
\end{lem}

\begin{proof}
Let $f \colon Y \to X$ be the canonical modification of $X$ (see \cite{Kol13}*{Theorem 1.31}), and put
\[
\Gamma := f^*(K_X + \Delta) - K_Y.
\]
Then $Y$ has canonical singularities, and in particular, it is $\mathbb{Q}$-Gorenstein, which proves (1). 
Since $K_X + \Delta$ is $\mathbb{R}$-Cartier, so is $\Gamma$.
We have $f_*\Gamma = \Delta \ge 0$, and $-\Gamma \equiv_f K_Y$ is $f$-ample.
By the negativity lemma, we have $\Gamma \ge 0$. 
Since $\Gamma \ge 0$ and $- \Gamma$ is $f$-ample, we have 
$\operatorname{Exc}(f) \subset \operatorname{Supp}(\Gamma)$.

We write 
\[
f^*S = S_Y + F
\]
with an effective exceptional divisor $F$. 
This $F$ has integer coefficients since $S$ is Cartier. 
Note that $(Y, f^* S + \Gamma)$ is log canonical since 
$(X, S + \Delta)$ is log canonical and $\Gamma \ge 0$. 
Therefore, it follows that $F = 0$ since 
$\operatorname{Supp} (F) \subset \operatorname{Exc}(f) \subset \operatorname{Supp}(\Gamma)$. 
Hence $f^*S = S_Y$. In particular, $S_Y$ is a Cartier prime divisor, $f^{-1}(S) = S_Y$, and (3) holds.

It remains to show that $S_Y$ is normal.
Since $Y$ has canonical singularities, it is Cohen--Macaulay, and hence so is the Cartier divisor $S_Y$.
In particular, $S_Y$ satisfies Serre's condition $S_2$.
By Serre's criterion for normality, it suffices to show that $S_Y$ satisfies $R_1$.
Let $\eta$ be a codimension-one point of $S_Y$.
If $\eta \notin \operatorname{Exc}(f)$, then $S_Y$ is regular at $\eta$ because $S$ is normal.
Suppose that $\eta \in \operatorname{Exc}(f)$.
Since $\eta \in \operatorname{Exc}(f) \subset \operatorname{Supp}(\Gamma)$ and $(Y, S_Y + \Gamma)$ is log canonical, we are in case (16.6.3) of \cite{92}*{16.6} (see the proof of \cite{92}*{16.7}), and $S_Y$ is normal at $\eta$. This completes the proof.
\end{proof}

\begin{prop}\label{prop:reduction_QGor}
Let $n$ be a positive integer. Suppose that Conjecture~\ref{conj:PIA} holds whenever $\dim X = n$ and $X$ is $\mathbb{Q}$-Gorenstein. Then Conjecture~\ref{conj:PIA} holds whenever $\dim X = n$.
\end{prop}

\begin{proof}
Let $(X, \Delta)$, $S$ and $x$ be as in Conjecture~\ref{conj:PIA} with $\dim X = n$. 
By Remark \ref{rmk:easy}, we may assume, after shrinking $X$ around $x$, that $(X, S + \Delta)$ is log canonical.

Let $f \colon Y \to X$ and $\Gamma$ be as in Lemma~\ref{lem:can_mod}, let $S_Y$ be the strict transform of $S$, let $f|_{S_Y} \colon S_Y \to S$ be the restriction, and set $W := f^{-1}(x)$.
By Lemma~\ref{lem:can_mod}~(1)–(3), $Y$ is $\mathbb{Q}$-Gorenstein of dimension $n$, $S_Y$ is a normal Cartier prime divisor, and  
\[
K_Y + S_Y + \Gamma = f^*(K_X + S + \Delta). 
\]
Note that $S_Y \not\subset \operatorname{Supp}(\Gamma)$ by the log canonicity of $(Y, S_Y + \Gamma)$. 
Define $\Gamma_{S_Y}$ by adjunction
\[
K_{S_Y} + \Gamma_{S_Y} = (K_Y + S_Y + \Gamma)|_{S_Y}.
\]
Then we have
\[
K_{S_Y} + \Gamma_{S_Y} = (f|_{S_Y})^*(K_S + \Delta |_S).
\]
Moreover, $W \subset S_Y$ and $W = (f|_{S_Y})^{-1}(x)$ by Lemma~\ref{lem:can_mod}~(2).

Since we have 
\[
\operatorname{mld}_x(X, S + \Delta) = \operatorname{mld}_W(Y, S_Y + \Gamma), \quad
\operatorname{mld}_x(S, \Delta |_S) = \operatorname{mld}_W(S_Y, \Gamma_{S_Y}), 
\]
it suffices to show
\[
\operatorname{mld}_{\eta_C}(Y, S_Y + \Gamma) = \operatorname{mld}_{\eta_C}(S_Y, \Gamma_{S_Y})
\]
for every irreducible closed subset $C \subset W$, where $\eta_C$ denotes the generic point of $C$.

Fix an irreducible closed subset $C \subset W$. 
Then by \cite{Amb99}*{Proposition 2.1}, there exists a closed point $y \in C$ such that
\begin{align*}
\operatorname{mld}_y(Y, S_Y + \Gamma) &= \operatorname{mld}_{\eta_C}(Y, S_Y + \Gamma) + \dim C, \\ 
\operatorname{mld}_y(S_Y, \Gamma_{S_Y}) &= \operatorname{mld}_{\eta_C}(S_Y, \Gamma_{S_Y}) + \dim C.
\end{align*}
Therefore, it suffices to show
\[
\operatorname{mld}_y(Y, S_Y + \Gamma) = \operatorname{mld}_y(S_Y, \Gamma_{S_Y}).
\]
Since $Y$ is $\mathbb{Q}$-Gorenstein of dimension $n$, this follows from the hypothesis that Conjecture~\ref{conj:PIA} holds for $\mathbb{Q}$-Gorenstein varieties of dimension $n$.
\end{proof}

\begin{rmk}\label{rmk:key_inclusion}
In the proof above, the inclusion $f^{-1}(x) \subset S_Y$ is essential: it allows us to compute both $\operatorname{mld}_x(X, S + \Delta)$ and $\operatorname{mld}_x(S, \Delta |_S)$ on the same closed subset $W = f^{-1}(x)$.
This inclusion is guaranteed by Lemma~\ref{lem:can_mod}~(2). 
If we instead take $f$ to be a dlt blow-up of $(X, S + \Delta)$, the variety $Y$ is still $\mathbb{Q}$-Gorenstein, but this inclusion is not guaranteed in general.
\end{rmk}

\section{Precise inversion of adjunction for quotients of local complete intersections}\label{section:quotient}

In this section, we prove Theorem~\ref{thm:quotient_PIA}, which states that Conjecture~\ref{conj:PIA} holds for finite quotients of local complete intersection singularities, provided that the Cartier divisor $S$ is klt.
We reduce it to the case where $G$ is cyclic by \cite{NS22}*{Remark 4.13}, and then to \cite{NS26a}*{Theorem 7.3} by constructing a suitable $G$-equivariant embedding (Lemma~\ref{lem:equiv_emb}).

We first prepare a lemma on finite covers, which will be used in the proofs of Theorem~\ref{thm:quotient_PIA} and Proposition~\ref{prop:PIA_QGor_dim4}. 

\begin{lem}\label{lem:cover_plt}
Let $\pi \colon X' \to X$ be a finite surjective morphism of normal varieties that is \'etale in codimension one.
Let $S$ be a normal prime Cartier divisor on $X$, and let $x' \in X'$ be a closed point such that $x := \pi(x') \in S$ and $\pi^{-1}(x) = \{x'\}$.
Suppose that $(X, S)$ is plt at $x$.
Then, after shrinking $X$ around $x$, the following hold:
\begin{enumerate}
\item $S' := \pi^*S$ is a normal prime divisor,
\item $(X', S')$ is plt, and
\item $\pi|_{S'} \colon S' \to S$ is finite, surjective and \'etale in codimension one, and
\[
K_{S'} = (K_{X'} + S')|_{S'} = (\pi|_{S'})^* K_S.
\]
\end{enumerate}
\end{lem}

\begin{proof}
Since $\pi$ is \'etale in codimension one, we have
\[
K_{X'} + \pi^*S = \pi^*(K_X + S).
\]
Since $(X, S)$ is plt at $x$, the pair $(X', \pi^*S)$ is also plt at $x'$ (see \cite{Kol13}*{Corollary 2.43}).
Shrinking $X$ around $x$, we may assume that $(X', \pi^*S)$ is plt.
Then $\pi^*S$ is normal by \cite{KM98}*{Proposition 5.51}, and hence, after further shrinking $X$ around $x$, it is a normal prime divisor.
This proves (1) and (2).

Since $S$ and $S'$ are normal Cartier divisors, we have
\[
K_{S'} = (K_{X'} + S')|_{S'} = \bigl(\pi^*(K_X + S)\bigr)|_{S'}
= (\pi|_{S'})^*\bigl((K_X + S)|_S\bigr) = (\pi|_{S'})^* K_S.
\]
Therefore $\pi|_{S'}$ is \'etale in codimension one, which proves (3).
\end{proof}

The following lemma shows that, for a cyclic group action, a local complete intersection singularity at a fixed point can be embedded equivariantly into an affine space with a linear action so that it is defined by semi-invariant equations. 

\begin{lem}\label{lem:equiv_emb}
Let $X'$ be a normal quasi-projective variety, let $G$ be a finite cyclic group acting on $X'$ freely in codimension one, and let $x' \in X'$ be a closed point fixed by $G$ such that $X'$ is a local complete intersection at $x'$.
Then, after replacing $X'$ by a $G$-stable affine open neighborhood of $x'$, there exist an affine space $W$ with a linear $G$-action and a $G$-equivariant closed immersion $X' \hookrightarrow W$ sending $x'$ to the origin such that
\begin{enumerate}
\item the defining ideal of $X'$ in $k[W]$ is generated by a $k[W]$-regular sequence of $G$-semi-invariant polynomials vanishing at the origin, and
\item the action of $G$ on $W$ is free in codimension one.
\end{enumerate}
\end{lem}

\begin{proof}
Since $X'$ is quasi-projective and $x'$ is fixed by $G$, the intersection of the $G$-translates of an affine open neighborhood of $x'$ is a $G$-stable affine open neighborhood of $x'$.
Replacing $X'$ by it, we may assume that $X' = \operatorname{Spec} B$ is affine.

We next present $B$ as a $G$-equivariant quotient of a polynomial ring.
Let $\mathfrak{m}_{x'} \subset B$ be the maximal ideal of $x'$. 
It is $G$-stable since $x'$ is fixed by $G$.
Choose elements $b_1, \dots, b_N \in \mathfrak{m}_{x'}$ generating $B$ as a $k$-algebra.
Let
\[
A := k[\, y_{g,i} \mid g \in G,\ 1 \le i \le N \,]
\]
be the polynomial ring in $N \cdot \# G$ variables, on which $G$ linearly acts by $\tau(y_{g,i}) := y_{\tau g,i}$ for $\tau \in G$.
Define a surjective $k$-algebra homomorphism $\varphi \colon A \to B$ by $\varphi(y_{g,i}) := g(b_i)$. 
By definition, $\varphi$ is $G$-equivariant. 

Let $\mathfrak{m} \subset A$ be the maximal ideal generated by all the variables $y_{g,i}$. 
Since $\varphi(\mathfrak{m}) \subset \mathfrak{m}_{x'}$, we have $\varphi^{-1}(\mathfrak{m}_{x'}) = \mathfrak{m}$.
Let $I := \ker \varphi$. Both $I$ and $\mathfrak{m}$ are $G$-stable.
Since $A_{\mathfrak{m}}/I_{\mathfrak{m}} \cong B_{\mathfrak{m}_{x'}}$ is a complete intersection and $A_{\mathfrak{m}}$ is regular, the ideal $I_{\mathfrak{m}}$ is generated by a regular sequence of length
\[
c := \operatorname{ht}(I_{\mathfrak{m}}).
\]
In particular, the minimal number of generators of $I_{\mathfrak{m}}$ is $c$, and hence
\[
C := I/\mathfrak{m}I \cong I_{\mathfrak{m}}/\mathfrak{m}I_{\mathfrak{m}}
\]
is a $c$-dimensional $k$-vector space. 
Since $I$ and $\mathfrak{m}$ are $G$-stable, $G$ acts linearly on $C$.

We now construct semi-invariant generators.
Let $r := \# G$, let $\sigma$ be a generator of $G$, and let $\zeta \in k$ be a primitive $r$-th root of unity.
Since $\sigma^r = 1$, the $\sigma$-action on $C$ is diagonalizable.
Thus we can choose a basis $\bar{f}_1, \dots, \bar{f}_c$ of $C$ and integers $0 \le a_i < r$ such that $\sigma(\bar{f}_i) = \zeta^{a_i} \bar{f}_i$ for $1 \le i \le c$. 
First, we choose any lift $u_i \in I$ of $\bar{f}_i$, and set
\[
f_i := \frac{1}{r} \sum_{j=0}^{r-1} \zeta^{-a_i j} \sigma^j(u_i) \in I.
\]
Since the image of $\sigma^j(u_i)$ in $C$ is $\sigma^j(\bar{f}_i) = \zeta^{a_i j} \bar{f}_i$, the image of $f_i$ in $C$ is
\[
\frac{1}{r} \sum_{j=0}^{r-1} \zeta^{-a_i j} \sigma^j(\bar{f}_i) = \frac{1}{r} \sum_{j=0}^{r-1} \bar{f}_i = \bar{f}_i.
\]
Thus $f_i$ is again a lift of $\bar{f}_i$. 
Moreover, each $f_i$ is $G$-semi-invariant since
\[
\sigma(f_i) = \frac{1}{r} \sum_{j=0}^{r-1} \zeta^{-a_i j} \sigma^{j+1}(u_i) = \zeta^{a_i} f_i. 
\]
By Nakayama's lemma, $f_1, \dots, f_c$ generate $I_{\mathfrak{m}}$.
Since $A_{\mathfrak{m}}$ is regular and $\operatorname{ht} I_{\mathfrak{m}} = c$, the elements $f_1, \dots, f_c$ form an $A_{\mathfrak{m}}$-regular sequence. 

We next pass from the local ring to a $G$-stable affine open neighborhood.
Since $I_{\mathfrak{m}} = (f_1, \dots, f_c) A_{\mathfrak{m}}$ and $f_1, \dots, f_c$ is an $A_{\mathfrak{m}}$-regular sequence, there exists $s \in A \setminus \mathfrak{m}$ such that $I A_s = (f_1, \dots, f_c) A_s$ and $f_1, \dots, f_c$ is an $A_s$-regular sequence.
Set $h := \prod_{g \in G} g(s)$.
Then $h$ is $G$-invariant, and $h \in A \setminus \mathfrak{m}$.
Since $A_h$ is a localization of $A_s$, we have $I A_h = (f_1, \dots, f_c) A_h$, and $f_1, \dots, f_c$ is an $A_h$-regular sequence.
Replacing $X'$ by the $G$-stable affine open neighborhood $\{\varphi(h) \neq 0\}$ of $x'$, that is, $B$ by $B_{\varphi(h)}$, we may assume that
\[
B \cong A_h/(f_1, \dots, f_c).
\]

We now rewrite $B$ as a quotient of a polynomial ring by a regular sequence.
Multiplying $h$ by a nonzero constant, we may assume $h - 1 \in \mathfrak{m}$.
Let $z$ be a new variable, and extend the action of $G$ to $A[z]$ by letting $G$ act trivially on $z$.
Then $H := (1 + z)h - 1$ is $G$-invariant, and we have $G$-equivariant isomorphisms
\[
B \cong A_h/(f_1, \dots, f_c) \cong A[z]/(H, f_1, \dots, f_c).
\]
Here, $H, f_1, \dots, f_c$ is an $A[z]$-regular sequence of $G$-semi-invariant polynomials
since $H$ is a non-zero-divisor in $A[z]$ and we have $A[z]/(H) \cong A_h$. 
Under the above isomorphism, the maximal ideal of $x'$ corresponds to $\mathfrak{m} A[z] + (z)$.

Set $W := \operatorname{Spec} A[z]$. Since $G$ permutes the variables $y_{g,i}$ and fixes $z$, this is an affine space with a linear $G$-action, and the above isomorphism gives a $G$-equivariant closed immersion $X' \hookrightarrow W$ sending $x'$ to the origin, whose defining ideal is generated by the regular sequence $H, f_1, \dots, f_c$ of $G$-semi-invariant polynomials. This proves (1).

Finally, we prove (2). Let $1 \neq g \in G$.
The fixed locus $W^g$ is a linear subspace of $W$ containing the origin. 
If $W^g \subset W$ is a hyperplane, then $\operatorname{codim}_{X'} (X' \cap W^g) \le 1$, 
which contradicts the assumption that the $G$-action on $X'$ is free in codimension one. 
Hence, the $G$-action on $W$ is free in codimension one.
\end{proof}

Using Lemmas~\ref{lem:cover_plt} and \ref{lem:equiv_emb}, we now prove the main result of this section.

\begin{thm}\label{thm:quotient_PIA}
Let $X'$ be a normal quasi-projective variety, let $G$ be a finite group acting on $X'$ freely in codimension one, and let $x' \in X'$ be a closed point such that $X'$ is a local complete intersection at $x'$.
Let $q \colon X' \to X := X'/G$ be the quotient morphism.
Let $(X, \Delta)$ and $S$ be as in Conjecture~\ref{conj:PIA} with $x := q(x') \in S$, and suppose that $S$ is klt at $x$.
Then Conjecture~\ref{conj:PIA} holds for $(X, \Delta)$, $S$ and $x$.
\end{thm}

\begin{proof}
Let $H \subset G$ be the stabilizer of $x'$, and let $x_H \in X'/H$ be the image of $x'$.
The induced morphism $p \colon X'/H \to X$ is \'etale at $x_H$, and $H$ acts on $X'$ freely in codimension one. 
Replacing $X$ by $X'/H$, we may assume that $G$ fixes $x'$. 
In particular, $q^{-1}(x) = \{x'\}$.

Note that $X'$ is Gorenstein at $x'$, since $X'$ is a local complete intersection at $x'$.
Since $q$ is \'etale in codimension one, it follows that $X$ is $\mathbb{Q}$-Gorenstein at $x$ (see \cite{KM98}*{Proposition 5.20(1)}).
Hence, after shrinking $X$ around $x$, we may assume that $X$ is $\mathbb{Q}$-Gorenstein.
Since $S$ is klt at $x$, the pair $(X, S)$ is plt at $x$ by inversion of adjunction (see Theorem~\ref{thm:IA}~(1)).
By Lemma~\ref{lem:cover_plt}, after shrinking $X$ around $x$, the divisor $S' := q^*S$ is a normal prime Cartier divisor on $X'$ and $q|_{S'} \colon S' \to S$ is \'etale in codimension one. 

Since $X$ is $\mathbb{Q}$-Gorenstein, we can apply Remark~\ref{rmk:ideal_setting}.
Let $\mathfrak{a}_\Delta$ be the $\mathbb{R}$-ideal sheaf on $X$ given there.
By Remark~\ref{rmk:ideal_setting}, it suffices to prove
\[\tag{$\diamondsuit$}
\operatorname{mld}_x\bigl(X, \mathcal{O}_X(-S)\mathfrak{a}_\Delta\bigr) = \operatorname{mld}_x(S, \mathfrak{a}_\Delta \mathcal{O}_S).
\]

We first reduce to the case where $G$ is cyclic.
For $\gamma \in G$, set $X_\gamma := X'/\langle \gamma \rangle$ and $S_\gamma := S'/\langle \gamma \rangle$, let $p_\gamma \colon X_\gamma \to X$ be the induced morphism, and let $x_\gamma \in S_\gamma$ be the image of $x'$.
Then $S_\gamma = p_\gamma^* S$ is a normal prime Cartier divisor on $X_\gamma$.
By \cite{NS22}*{Remark 4.13}, applied to the $G$-actions on $X'$ and $S'$, we have
\begin{align*}
\operatorname{mld}_x\bigl(X, \mathcal{O}_X(-S)\mathfrak{a}_\Delta\bigr) &= \min_{\gamma \in G} \operatorname{mld}_{x_\gamma}\bigl(X_\gamma, \mathcal{O}_{X_\gamma}(-S_\gamma)\,\mathfrak{a}_\Delta \mathcal{O}_{X_\gamma}\bigr), \\
\operatorname{mld}_x(S, \mathfrak{a}_\Delta \mathcal{O}_S) &= \min_{\gamma \in G} \operatorname{mld}_{x_\gamma}\bigl(S_\gamma, \mathfrak{a}_\Delta \mathcal{O}_{S_\gamma}\bigr).
\end{align*}
Hence it suffices to prove
\[
\operatorname{mld}_{x_\gamma}\bigl(X_\gamma, \mathcal{O}_{X_\gamma}(-S_\gamma)\,\mathfrak{a}_\Delta \mathcal{O}_{X_\gamma}\bigr) = \operatorname{mld}_{x_\gamma}\bigl(S_\gamma, \mathfrak{a}_\Delta \mathcal{O}_{S_\gamma}\bigr)
\]
for each $\gamma \in G$. 
Since $\langle \gamma \rangle$ acts on $X'$ freely in codimension one and fixes $x'$, and since $S_\gamma$ is a normal prime Cartier divisor on $X_\gamma$ which is klt at $x_\gamma$ by \cite{Kol13}*{Corollary 2.43}, we may replace $G$ by $\langle \gamma \rangle$ and assume that $G$ is cyclic.

By Lemma~\ref{lem:equiv_emb}, after a $G$-stable shrinking of $X'$ around $x'$, there exists a $G$-equivariant closed immersion $X' \hookrightarrow W$ into an affine space $W$ with a linear $G$-action, sending $x'$ to the origin, such that the defining ideal of $X'$ is generated by a regular sequence of $G$-semi-invariant polynomials and the action of $G$ on $W$ is free in codimension one.
We may therefore apply \cite{NS26a}*{Theorem 7.3} to $X = X'/G$, the klt Cartier divisor $S$ and the $\mathbb{R}$-ideal sheaf $\mathfrak{a}_\Delta$, which gives ($\diamondsuit$).
This completes the proof.
\end{proof}

\begin{rmk}\label{rmk:comparison_NS}
Theorem~\ref{thm:quotient_PIA} improves \cite{NS26a}*{Theorem 7.3} in the following two points.
In \cite{NS26a}*{Theorem 7.3}, it is assumed that
\begin{itemize}
\item $X'$ is a closed subvariety of an affine space $\mathbb{A}^N_k$ with a linear $G$-action, defined by $c := N - \dim X'$ $G$-semi-invariant equations, and
\item $G$ acts freely in codimension one also on $S' = q^*S$.
\end{itemize}
The first assumption is replaced by the condition that $X'$ is a local complete intersection at $x'$: by Lemma~\ref{lem:equiv_emb}, the ambient space need not be an affine space from the beginning.
The second assumption follows from the fact that $S$ and $S'$ are normal Cartier divisors (see Lemma~\ref{lem:cover_plt}).
\end{rmk}

\section{Precise inversion of adjunction in dimension four}\label{section:dim4}

In this section, we treat the case where $\operatorname{mld}_x(X, S + \Delta) \ge 1$ in dimension four, and complete the proof of Theorem~\ref{thm:PIA_dim4}, the main theorem of this paper. 

\begin{prop}\label{prop:PIA_QGor_dim4}
Conjecture~\ref{conj:PIA} holds when $\dim X = 4$, $X$ is $\mathbb{Q}$-Gorenstein, and $\operatorname{mld}_x(X, S + \Delta) \ge 1$.
\end{prop}

\begin{proof}
Let $(X, \Delta)$, $S$ and $x$ be as in Conjecture~\ref{conj:PIA} with $\dim X = 4$, and suppose that $X$ is $\mathbb{Q}$-Gorenstein and $\operatorname{mld}_x(X, S + \Delta) \ge 1$.
By Remark~\ref{rmk:easy}, we have
\[
\operatorname{mld}_x(S, \Delta |_S) \ge \operatorname{mld}_x(X, S + \Delta) \ge 1,
\]
and the assertion is clear if $\operatorname{mld}_x(S, \Delta |_S) = 1$.
We may therefore assume $\operatorname{mld}_x(S, \Delta |_S) > 1$.
Note that $S$ is $\mathbb{Q}$-Gorenstein, since $K_S = (K_X + S)|_S$ and $X$ is $\mathbb{Q}$-Gorenstein.
Since $\Delta |_S$ is effective, we have
\[\tag{$\spadesuit$}
\operatorname{mld}_x(S) \ge \operatorname{mld}_x(S, \Delta |_S) > 1.
\]

We first show that $S$ is klt at $x$.
By ($\spadesuit$), $S$ is log canonical at $x$.
Suppose that $S$ is not klt at $x$, and let $W \subset S$ be a log canonical center of $S$ containing $x$.
Since $\dim S = 3$, $W$ is either the point $x$ or a curve.
In the former case, we have $\operatorname{mld}_x(S) = 0$. 
In the latter case, we have 
\[
\operatorname{mld}_x(S) \le \operatorname{mld}_{\eta_W}(S) + 1 = 1
\]
by the lower semi-continuity of minimal log discrepancies in dimension $3$ (\cite{Amb99}*{Theorem 0.1}). 
Both contradict ($\spadesuit$), so $S$ is klt at $x$.
Furthermore, $(X, S)$ is plt at $x$ by inversion of adjunction (see Theorem~\ref{thm:IA}~(1)).

Let $r$ be the Cartier index of $K_X$ at $x$.
After shrinking $X$ to an affine neighborhood of $x$, we may assume that $\mathcal{O}_X(rK_X) \cong \mathcal{O}_X$.
Let $q \colon X' \to X$ be the associated index-one cover \cite{KM98}*{Definition 5.19}.
Then $X'$ is a normal affine variety, $K_{X'}$ is Cartier, the cyclic group $G := \mu_r$ acts on $X'$ with $X = X'/G$, and $q$ is \'etale in codimension one.
Since $r$ is the Cartier index at $x$, the fiber $q^{-1}(x)$ consists of a single point $x'$.
By Lemma~\ref{lem:cover_plt}, after shrinking $X$ around $x$, the divisor $S' := q^*S$ is a normal prime Cartier divisor on $X'$, and
\[
K_{S'} = (K_{X'} + S')|_{S'} = (q|_{S'})^* K_S.
\]
In particular, $K_{S'}$ is Cartier.

We next show that $X'$ is smooth or has a hypersurface singularity at $x'$.
Since $K_{S'} = (q|_{S'})^*K_S$, \cite{Kol13}*{2.42} and ($\spadesuit$) give
\[
\operatorname{mld}_{x'}(S') \ge \operatorname{mld}_x(S) > 1.
\]
Since $K_{S'}$ is Cartier, $\operatorname{mld}_{x'}(S')$ is an integer,
and hence $S'$ has a Gorenstein canonical singularity at $x'$.
By the classification of three-dimensional Gorenstein canonical singularities (see \cite{Rei80}*{Theorem 2.2}; see also \cite{Amb99}*{Lemma 3.2 and Proposition 3.3}), $S'$ is either smooth or has a cDV singularity at $x'$.
In particular, $\widehat{\mathcal{O}}_{S', x'}$ is regular or a hypersurface ring.
Let $R := \widehat{\mathcal{O}}_{X', x'}$, and let $f \in R$ be a local equation of $S'$ at $x'$. Then $f$ is a non-zero-divisor in $R$, and $R/(f) \cong \widehat{\mathcal{O}}_{S', x'}$. 
Therefore, $R$ is a complete intersection (see \cite{BH93}*{Theorem 2.3.4}).
Moreover, denoting by $\operatorname{edim}$ the embedding dimension, we have
\[
\operatorname{edim} R - \dim R \le \bigl(\operatorname{edim} R/(f) + 1\bigr) - \bigl(\dim R/(f) + 1\bigr) \le 1, 
\]
which proves that $X'$ is smooth or has a hypersurface singularity at $x'$.

Therefore the hypotheses of Theorem~\ref{thm:quotient_PIA} are satisfied by $X'$, $G$, $x'$ and $q$, and $S$ is klt at $x$.
Theorem~\ref{thm:quotient_PIA} gives the assertion.
\end{proof}

Combining the results so far, we obtain the main theorem of this paper.

\begin{thm}\label{thm:PIA_dim4}
Conjecture~\ref{conj:PIA} holds when $\dim X = 4$.
\end{thm}

\begin{proof}
By Proposition~\ref{prop:reduction_QGor} applied with $n = 4$, it suffices to prove Conjecture~\ref{conj:PIA} when $\dim X = 4$ and $X$ is $\mathbb{Q}$-Gorenstein.
If $\operatorname{mld}_x(X, S + \Delta) < 1$, this follows from Proposition~\ref{prop:PIA_low}.
Otherwise, it follows from Proposition~\ref{prop:PIA_QGor_dim4}.
\end{proof}

As a byproduct of the proof of Proposition~\ref{prop:PIA_QGor_dim4}, we obtain the following structure theorem for four-dimensional singularities admitting a normal Cartier divisor with minimal log discrepancy greater than one.

\begin{thm}\label{thm:hyperquotient}
Let $X$ be a normal $\mathbb{Q}$-Gorenstein variety with $\dim X = 4$, let $S$ be a normal prime Cartier divisor on $X$, and let $x \in S$ be a closed point.
Suppose that $\operatorname{mld}_x(X, S) > 1$.
Then $X$ has a hyperquotient singularity at $x$. More precisely, the index-one cover of $X$ at $x$ is smooth or has a hypersurface singularity at the point over $x$.
\end{thm}

\begin{proof}
By Remark~\ref{rmk:easy}, we have $\operatorname{mld}_x(S) \ge \operatorname{mld}_x(X, S) > 1$, and the assertion follows from the proof of Proposition~\ref{prop:PIA_QGor_dim4} with $\Delta = 0$.
\end{proof}

\end{document}